
\input amstex
\loadbold


\define\jsp{J.\ Statist.\ Phys.}
\define\aap{Ann.\ Appl.\ Probab.}
\define\ap{Ann.\ Probab.}

\define\zwvg{Z.\ Wahrsch.\ Verw.\ Gebiete}

\define\<{\langle}
\define\>{\rangle}

\define\({\left(}
\define\){\right)}
\define\[{\left[}                                   
\define\]{\right]}                                            
\define\lbrak{\bigl\{}
\define\rbrak{\bigr\}}
\define\lbrakk{\biggl\{}
\define\rbrakk{\biggr\}}

\define\rate{r}  



\define\kappatil{{\tilde \kappa}}
\define\latil{{\tilde \lambda}}

\define\vtil{{\tilde v}}
\define\wtil{{\widetilde w}}



\define\e{\varepsilon}

\define\ind{{\boldkey 1}}

\define\epi#1{\text{epi}\,#1}

\define\pikkuhyppy{\vskip .1in}


\define\mmN{{\bold N}}

\define\mmR{{\bold R}}

\define\mmZ{{\bold Z}}

\define\mmu{{\bold u}}
\define\mmv{{\bold v}}

\define\mmx{{\bold x}}
\define\mmy{{\bold y}}


\define\XX{{\Cal X}}

\define\DD{{\Cal D}}

\define\LL{{\Cal L}}

\define\VV{{\Cal V}}


\documentstyle{amsppt}
\magnification=\magstep1
\document     
\baselineskip=12pt
\pageheight{43pc} 

\NoBlackBoxes

\centerline{\bf Hydrodynamic profiles for the totally asymmetric
exclusion process}

\pikkuhyppy

\centerline{\bf with a slow bond}

\pikkuhyppy

$$\text{2000}$$

\hbox{}

\centerline{  Timo  Sepp\"al\"ainen
\footnote""{ Research partially supported by NSF grant DMS-9801085. }}
\hbox{}
\centerline{Department of Mathematics}
\centerline{Iowa State University}
\centerline{Ames, Iowa 50011, USA}
\centerline{seppalai\@iastate.edu}

\hbox{}

\hbox{}

\flushpar
{\it Summary.} We study a totally asymmetric simple 
exclusion process where jumps happen at rate one,
except at the origin where the rate is lower.
We prove a hydrodynamic scaling limit to a macroscopic
profile described by a variational formula. The limit
is valid for all values of the slow rate. The 
only assumption required is that a law
of large numbers holds for
the initial particle distribution. This includes
also deterministic initial configurations. 
The hydrodynamic description contains as an 
unknown parameter the macroscopic rate at the 
origin, which is strictly larger than the 
microscopic slow rate. The limit is proved by
the variational coupling method. 


\hfil

\hfil


\hfil

\flushpar

\hfil

\flushpar
Short Title: Exclusion with a slow bond 

\break

\flushpar
{\bf 1. Introduction}

\hbox{}

The  exclusion process 
 is one of the most basic  stochastic
interacting particle  processes.
It consists of particles executing independent 
random walks on a graph, subject to the 
exclusion interaction that prevents two particles
from occupying the same vertex. 
In this paper we look at the version where 
the underlying graph is the integer lattice 
$\mmZ$ with nearest-neighbor bonds. The  particles 
take only nearest-neighbor steps to the right, 
and a jump is permitted only if the next site to
the right is vacant. In the  spatially homogeneous 
case jumps across all nearest-neighbor bonds happen
at exponential rate one. This model goes by the name 
TASEP, or totally asymmetric simple exclusion process. 

In the early 1990's Janowsky and Lebowitz \cite{JL1, JL2}
introduced an interesting variant of TASEP. 
Jumps across a specific bond, say $\{0,1\}$, happen
at a rate $r$ strictly lower than the rate 1 everywhere
else. Imagine for example a tollbooth on a 
single lane highway that slows down the flow of cars. 
For the slow bond model one investigates the same 
questions 
 as for homogeneous TASEP, such as large scale behavior 
and invariant distributions. 
Simulations and partial results
are available, but crucial results remain open.

In the present paper we prove a hydrodynamic limit 
for the slow bond model. In one sense our 
result is complete and general. It is valid for all
values $r\in(0,1)$ of the slow rate. We need no 
extra assumptions on the initial distribution of the 
particles, only that the 
law of large numbers for the density profile
 be valid at time zero. But the
hydrodynamic limit cannot make the most
interesting distinction, namely whether the slow
bond disturbs hydrodynamic profiles for all values
$r<1$. This question can be resolved only through 
sharper control of particle-level evolution. 

The  macroscopic effect
of the slow bond can be described like this:
Corresponding to the microscopic rate profile
(rate identically 1 at all sites except $0$,
and rate equal to $r$ at $0$) there is a macroscopic
rate profile $\lambda(x)$ defined for $x\in\mmR$
 such that $\lambda(x)=1$
for $x\ne 0$, and 
$\lambda(0)\in(r,1]$. 
The macroscopic evolution of the particle density
is the solution of an optimal control problem whose
running cost depends on the space variable 
$x$ through the function $\lambda(x)$. 
If $\lambda(0)=1$, there is no 
visible macroscopic disturbance from the slow bond. 
The interesting open question becomes:
Is $\lambda(0)<1$ for all $r<1$? 

We give a precise definition of the quantity $\lambda(0)$ 
in terms of a last-passage growth model whose
association with TASEP is well-known. In this
formulation the open problem can be cast in a 
percolation-type framework, and perhaps made amenable to
 other techniques, such as the random matrix
methods that have recently been applied to growth models
with great success.

To prove the result we use the variational coupling
method initiated in 
\cite{Se1} which, when it applies, gives laws of large 
numbers without any knowledge of invariant distributions. 
This  approach allows us to work  directly
on the slow bond model, 
without intermediate approximations with more
tractable models. The proof extends to the situation
with multiple slow bonds, and also works for more
general exclusion processes that admit $K$ particles
per site. 

In addition to the Janowsky--Lebowitz papers mentioned
above, there is a handful of other related work. 
Covert and Rezakhanlou \cite{CR} derived a bound 
for the critical value of the slow rate by approximating
the slow bond model with an exclusion process whose
rates vary more regularly in space. Liggett's \cite{Li}
new monograph discusses the slow bond model and 
proves bounds for the critical value of the rate. 
As for Janowsky and Lebowitz, 
the approach is through finite systems with open 
boundaries, with system size tending to infinity. 
 The case of a zero-range
process (ZRP) with a slow site was treated by Landim \cite{La}. 
In ZRP the particles accumulate at the slow site
and produce a point mass in the hydrodynamic
profile. The ZRP is a more tractable model because it
retains product-form invariant distributions even when 
rates lose spatial homogeneity. 

Our paper is organized as follows. The main result 
is the hydrodynamic limit Thm.\ 2.2. On the way to 
it we define the macroscopic rate $\lambda(0)$ 
in terms of the growth model, and state some
bounds for it in Thm.\ 2.1.  As corollaries to
Thm.\ 2.2 we compute the macroscopic profiles that
evolve from constant initial profiles, characterize 
the macroscopically invariant profiles in the 
range $[\rho^*,1-\rho^*]$ between the lower and upper
critical densities, and prove one property  of 
invariant measures. These results are in Sect.\ 2.
The remainder of the paper contains the proofs.

\hbox{}

\flushpar
{\bf 2. Results}

\hbox{}

The process operates according to these rules: 
Indistinguishable particles occupy the sites of the 
one-dimensional integer lattice $\mmZ$. Each site
has at most
one particle, so the state of the process is described
by the  occupation numbers 
$\eta=(\eta_i)_{i\in\mmZ}$, where $\eta_i=1$ if site
$i$ is occupied, and $0$ if site $i$ is empty. Particles
take nearest-neighbor steps to the right, subject to the 
exclusion rule, at exponential rate 1, with one exception:
jumps from site $0$ to site $1$ happen at rate $\rate$.
The rate $\rate$ is a fixed constant in the range $(0,1]$. 
On the compact state space $\{0,1\}^\mmZ$, the dynamics
has the infinitesimal generator 
$$Lf(\eta)=\sum_{i\ne 0} \eta_i(1-\eta_{i+1})
[f(\eta^{i,i+1})-f(\eta)]\ +\ \rate\eta_0(1-\eta_{1})
[f(\eta^{0,1})-f(\eta)]\,,
\tag 2.1
$$
where $\eta^{i,i+1}=\eta-\delta_i+\delta_{i+1}$ denotes
the configuration that results after a single particle 
jumps from $i$ to $i+1$. When $r=1$ this is the generator
of TASEP, the totally asymmetric simple exclusion 
process. 

We prove a hydrodynamic limit for this process in the
usual Euler scale, 
for all values of $r$. As is well-known, when $r=1$ the 
 macroscopic particle density $\rho(x,t)$ 
obeys 
the scalar conservation 
law 
$$ \rho_t + f_0(\rho)_x=0
\tag 2.2
$$
with current 
$$
f_0(\rho)=\rho(1-\rho). 
\tag 2.3
$$
Starting with Rost \cite{Ro} in 1981, this hydrodynamic
limit of the space-homogeneous TASEP has gone 
through many stages of 
generalization and refinement. 
See Ch.\ 8 of \cite{KL}, Part III of \cite{Li}
 and their notes and references. 

The slow bond restricts
the range of admissible currents. In TASEP the range
of currents is $[0,1/4]$, and the maximal
current $1/4=f_0(1/2)$ occurs at density $\rho=1/2$.  
In the slow bond model the maximal current is some 
value $J_{\max}(\rate)\le 1/4$. We know rigorously that 
$J_{\max}(\rate)< 1/4$ if $\rate$ is small enough. When this
happens, the densities $\rho$ around $1/2$ for 
which $f_0(\rho)>J_{\max}(\rate)$ become inadmissible. 
This we can prove on the hydrodynamic scale.

Let 
$\kappa(\rate)=J_{\max}(\rate)^{-1}$ be the reciprocal 
of the maximal rate. Thus 
$\kappa(\rate)$ is a  nonincreasing 
 function of $\rate\in(0,1]$, with values in the range
 $[4,\infty)$ and $\kappa(1)=4$. 
Since we do not know 
the exact dependence of $\kappa(\rate)$ on  $\rate$, 
we shall give it a precise definition through 
 a well-known growth model.

\hbox{}

\flushpar
{\it 2.1. Definition of $\kappa(\rate)$}

\hbox{}

Consider the following last-passage growth model
on the first quadrant of the plane.  
 For $(i,j)\in\mmN^2$, let
$Y_{i,j}$ be i.i.d.\ exponentially distributed 
random variables with common expectation $E[Y_{i,j}]=1$. 
Let the deterministic weights $w_{i,j}$ be given by 
$$
w_{i,j}=\cases
1\,, &i\ne j\\
1/{\rate}\,, &i=j.
\endcases
\tag 2.4
$$
For $n=1,2,3,\ldots$, let the passage time of 
site $(n,n)$ be 
$$T_{n,n}=\max_\pi \sum_{k=1}^{2n-1}
w_{i_k,j_k}Y_{i_k,j_k}\,,
\tag 2.5
$$
where the maximum is over paths 
$\pi=\{(1,1)=(i_1,j_1),(i_2,j_2),\ldots, 
(i_{2n-1},j_{2n-1})=(n,n)\}$ that take 
steps only to the right and up:
for each $k$, 
$$\text{$(i_k,j_k)-(i_{k-1},j_{k-1})=(1,0)$ or $(0,1)$.}
\tag 2.6
$$
An obvious superadditivity holds for the random  variables 
$\{T_{n,n}\}$, and  one can derive moment bounds 
(see for example Thm.\ 6.3 in \cite{GW}, or Prop.\ 5.1 in \cite{Se3})
sufficient for Kingman's
\cite{Ki} subadditive ergodic 
theorem, to obtain a strong law of large numbers:
There exists a constant $\kappa(\rate)$ such that
$$\lim_{n\to\infty}\frac1n T_{n,n} =\kappa(\rate)\ \text{ a.s.}
\tag 2.7
$$
This limit is taken as the definition of $\kappa(\rate)$. 
It is known that $\kappa(1)=4$. This is a special case of the
full interface result:  for $(x,y)\in\mmR_+^2$ 
$$\lim_{n\to\infty}\frac1n T_{[nx],[ny]} =
\bigl(\sqrt{x}+\sqrt{y}\bigr)^2\ \text{ a.s.\ in the case
$\rate=1$.}
\tag 2.8
$$
Proofs of (2.8) can be found in \cite{Ro} and \cite{Se1}. 

Currently we have this information on 
$\kappa(\rate)$:

\proclaim{Theorem 2.1}
The function  $\kappa(\rate)$, $0<\rate\le 1$,
 is continuous and nonincreasing. It satisfies these
bounds:
$$0\le \kappa(\rate_1)-\kappa(\rate_2)\le \frac1{\rate_1}
-\frac1{\rate_2}\ \ \ \text{for $\rate_1<\rate_2$, }
$$
and 
$$\max\biggl\{ 4\,,\, \frac32+\frac{\rate^2+2(1+\rate)}
{2\rate(1+\rate)}\biggr\} \le \kappa(\rate)
\le 3+ \frac1\rate\,.
\tag 2.9
$$
\endproclaim

The proof shows that 
neither bound in (2.9) is optimal. Janowsky and 
Lebowitz \cite{JL2} give the bound $J_{\max}(r)\le r(1-r)$ 
for $r\le 1/2$. (For the proof, see
 p.\ 277 in \cite{Li}.) For 
$\kappa(\rate)$ this gives 
$$\kappa(\rate)\ge \frac1{\rate(1-\rate)}\ \ \ \text{for $r\le 1/2$.}
\tag 2.10
$$
This
improves (2.9) 
in the range $r\in \bigl[ (\sqrt{17}-1)/8, 1/2 \bigr]$. 
However, the definition of $J_{\max}(r)$ in
 \cite{JL2} and \cite{Li}
 is not the same as ours, which is tied to the 
hydrodynamic limit. 
In
 \cite{JL2, Li},
  $J_{\max}(r)$  is the limiting  stationary current 
of a finite system with open boundaries, as the size
of the system tends to infinity. Of course 
the two quantities ought to be the same. Once 
optimal bounds are found with some approach, it will
be of interest to verify that the different 
definitions are the same. 

 Let us 
define 
$$\rate^*=\inf\{\rate\in(0,1]: \kappa(\rate)=4 \}\,.
$$
The interesting open problem is whether $\rate^*=1$.
In other words, does the macroscopic passage time 
$\kappa(\rate)$ rise strictly above 
$\kappa(1)=4$ as soon as $\rate<1$? Equivalently, is there
a forbidden range of densities around $1/2$ 
 as soon as $\rate<1$? 
 Simulations 
by Janowsky and Lebowitz \cite{JL2} suggest that such is the 
case. 
The lower bound in (2.9) gives  $\rate^*\ge (\sqrt{41}-3)/8\approx 0.425$, 
while (2.10) gives $r^*\ge 1/2$. 

\hbox{}

\flushpar
{\it 2.2. The hydrodynamic limit}

\hbox{}

Let $\rho_0$ be a given measurable function on $\mmR$ 
such that $0\le \rho_0(x)\le 1$ for all $x\in\mmR$. 
This is the initial macroscopic profile. 
Let $v_0$ be the antiderivative of $\rho_0$ defined by
$$v_0(0)=0\,,\, v_0(b)-v_0(a)=\int_a^b \rho_0(x)dx\,.
\tag 2.11
$$

Define the ``macroscopic rate profile''
 $\lambda(x)$  as 
$$
\lambda(x)=\cases
1\,, &x\ne 0\\
\dsize\frac4{\kappa(\rate)}\,, &x=0, 
\endcases
\tag 2.12
$$
where $\kappa(\rate)$ is defined by (2.7). 
For 
$x\in\mmR$, let 
$$g_0(x)=\sup_{0\le \rho\le 1}\{f_0(\rho)-x\rho\}$$ 
denote the Legendre conjugate of $f_0$. It is given by 
$$g_0(x)=\cases
-x\,, &x\le -1\\
(1/4)(1-x)^2\,,&-1\le x\le 1\\
0\,, &x\ge 1\,.
\endcases
\tag 2.13
$$
For  $x\in\mmR$ and $t\ge 0$, define $v(x,0)=v_0(x)$, and 
for $t>0$, 
$$v(x,t)=
\sup_{w(\cdot)}\lbrakk 
v_0(w(0))-\int_0^t \lambda(w(s))g_0\biggl(\frac{w'(s)}{\lambda(w(s))}\biggr)
ds \rbrakk\,.
\tag 2.14
$$
The supremum is over piecewise $C^1$ paths 
$w:[0,t]\to\mmR$ that satisfy $w(t)=x$. 
In Section 5.1 we give a formula for the path $w(\cdot)$
that minimizes the integral part inside the braces
 for a given initial
point $w(0)=q$. 

For each fixed $t$, $v(\cdot\,,t)$ is a Lipschitz function. 
Its $x$-derivative 
$$\rho(x,t)=\frac{\partial}{\partial x} v(x,t)
\tag 2.15
$$ represents the 
macroscopic density profile of the particles. 

If $\lambda(x)= 1$ for all $x$,  (2.14)--(2.15) give the 
entropy solution of (2.2) with initial data $\rho_0$. 
 This is the
density profile of space-homogeneous TASEP. 
Formula (2.14) expresses the sense in which $\lambda(x)$ can 
be regarded as the macroscopic rate profile. From (2.9) 
 we get 
$$\lambda(0)\ge \frac{4\rate}{3\rate+1}.
$$
This   implies that $\lambda(0)>\rate$ for $\rate<1$.
In other words,  the microscopic averaging does {\it not} simply 
reproduce the rate $\rate$ at the 
macroscopic level.

For the hydrodynamic limit, assume that we  have constructed
a sequence of exclusion processes $\eta^n(t)=
(\eta^n_i(t):i\in\mmZ)$ for times $t\ge 0$, where
$n=1,2,3,\ldots$ is the index of the sequence.
The initial distributions of the processes are 
arbitrary, subject to the condition that this weak law
of large numbers is valid: 
$$
\aligned
&\text{for all $a<b$ and $\e>0$, }\\
&\lim_{n\to\infty}P^n\lbrakk \biggl|\,
\frac1n\sum_{i=[na]+1}^{[nb]}\eta^n_i(0)- \int_a^b \rho_0(x)dx
\,\biggr|\ge \e\rbrakk=0\,.
\endaligned
\tag 2.16
$$
We wrote $P^n$ for the probability measure 
on the probability space of the process $\eta^n$. 
Let $J^n_i(t)$ denote the number of particles that have 
made the jump from site $i$ to $i+1$ in the time 
interval $[0,t]$, in the process $\eta^n(\cdot)$.

\proclaim{Theorem 2.2} Let the slow rate $r$ be any 
number in $(0,1]$. Under assumption {\rm (2.16)},
these weak laws of large numbers hold at macroscopic
times $t>0$:  for all real numbers $a<b$ and $\e>0$, 
$$
\lim_{n\to\infty}P^n\lbrak \bigl|\,
n^{-1}J^n_{[na]}(nt)-\bigl(v_0(a)-v(a,t)\bigr) 
\,\bigr|\ge \e\rbrak=0\,,
\tag 2.17
$$
and
$$
\lim_{n\to\infty}P^n\lbrakk \biggl|\,
\frac1n\sum_{i=[na]+1}^{[nb]}\eta^n_i(nt)-\int_a^b \rho(x,t)dx
\,\biggr|\ge \e\rbrakk=0\,,
$$
where 
 $v(x,t)$ is defined by {\rm (2.14)} and $\rho(x,t)=v_x(x,t)$. 
\endproclaim

 Let us derive 
some corollaries of Theorem 2.2. 
When $\lambda(0)<1$ there is a critical density
$$\rho^*=\frac12-\frac12\sqrt{1-\lambda(0)\,}\,,
\tag 2.18
$$
determined by the current condition $f_0(\rho^*)=\lambda(0)/4$. 
The slow bond disturbs the hydrodynamic profile 
only for densities in the range $\rho\in(\rho^*,1-\rho^*)$. 
When $\lambda(0)=1$ [or, equivalently, 
the passage time $\kappa(\rate)=4$], $\rho^*=1/2$
and the interval $(\rho^*,1-\rho^*)$ is empty. 

As a first application of Thm.\ 2.2 we derive the
macroscopic evolution of constant initial
profiles. 

\proclaim{Corollary 2.1} Assume the initial
distribution of the particles is chosen so that {\rm (2.16)} is
satisfied for the constant profile
   $\rho_0(x)\equiv \rho$. Then the 
macroscopic limits in Thm.\ 2.2 are described as 
follows:

{\it Case 1:}  $\rho\in[0,\rho^*]\cup[1-\rho^*,1]$. 
Then $v(x,t)=\rho x-t\rho(1-\rho)$ and 
macroscopically we see a constant density
$\rho(x,t)\equiv \rho$, exactly as for the homogeneous
TASEP. 

{\it Case 2:} $\lambda(0)<1$ and $\rho\in(\rho^*,1-\rho^*)$. 
The behavior changes around the blockage at the origin:
$$
v(x,t)=
\cases (1-\rho^*)x-t\lambda(0)/4\,,& -t(\rho-\rho^*)\le x\le 0\\
\rho^*x-t\lambda(0)/4\,,&0<x\le t(1-\rho^*-\rho)\,.
\endcases
\tag 2.19
$$
 Case 1 is still
valid for $x$ outside the range in {\rm (2.19)}. 
Correspondingly, the density profile has constant segments of 
upper and lower critical densities around the origin:
$$
\rho(x,t)=
\cases \rho\,, &\text{$x<-t(\rho-\rho^*)$ or $x>t(1-\rho^*-\rho)$}\\
1-\rho^*\,,& -t(\rho-\rho^*)< x< 0\\
\rho^*\,,&0<x< t(1-\rho^*-\rho)\,.
\endcases
\tag 2.20
$$
\endproclaim

The reason for this behavior is that the
maximal current permitted by the blocked system  is 
$\lambda(0)/4=\rho^*(1-\rho^*)$, while the unblocked TASEP has current
$f_0(\rho)=\rho(1-\rho)$. The blockage does not
disturb the system unless the system tries to transport
particles at a current above $\lambda(0)/4$. This happens
if $\rho(1-\rho)>\lambda(0)/4$, which is 
equivalent to $\rho\in(\rho^*,1-\rho^*)$. 

\proclaim{Corollary 2.2} Suppose $\lambda(0)<1$ so that 
$\rho^*\in(0,1/2)$. Let  $\rho_0$ be a  macroscopic profile
that is piecewise continuous in each bounded interval, and 
satisfies $\rho^*\le \rho_0(x)\le 1-\rho^*$. Suppose
$\rho_0$ is 
invariant under the macroscopic dynamics,
so that $\rho(x,t)=\rho_0(x)$
for a.e.\ $x$, for all times $t$. Then 
Lebesgue-almost everywhere 
$\rho_0(x)$ $\in$ $\{\rho^*, 1-\rho^*\}$, and $\rho_0$
takes 0--3 jumps according to these restrictions:
a jump from $1-\rho^*$ to $\rho^*$ can occur only 
at the origin, while  jumps from $\rho^*$ to $1-\rho^*$ 
can occur at any $x\in\mmR$.
\endproclaim

 We can control $\rho_0$ only almost everywhere because
a macroscopic density profile
$\rho(\cdot\,,t)$ is determined only through its integral
 $v(\cdot\,,t)$. The upward jumps
from $\rho^*$ to $1-\rho^*$ are the usual entropy shocks
of equation (2.2). The slow bond permits a downward
jump from $1-\rho^*$ to $\rho^*$ 
at the origin, which  violates the entropy condition
for the concave current $f_0$. 

By Cor.\ 2.2 there cannot be a
macroscopically 
invariant profile strictly in the range $(\rho^*,1-\rho^*)$.
Consequently, the process  cannot have  an invariant 
probability distribution 
for which a macroscopic profile exists in the sense of
(2.16) and lies in the range $(\rho^*,1-\rho^*)$.
The  interesting open problem is to prove the existence
of an invariant probability measure that corresponds to the 
macroscopically invariant non-entropy
 shock profile 
$\rho_0(x) =(1-\rho^*)\ind\{x<0\}+\rho^*\ind\{x>0\}$. 
About invariant measures we have  this to
say:

\proclaim{Corollary 2.3} For any value
$r$ of the slow rate, and for any 
$\rho\in[0,\rho^*]\cup[1-\rho^*,1]$, there exists
an invariant distribution $\mu$ for the process such that
$$\mu\{\eta_i=1,\eta_{i+1}=0\}=
\cases r^{-1}\rho(1-\rho)\,, &i= 0\\
\rho(1-\rho)\,, &i\ne 0.
\endcases
\tag 2.21
$$ 

In particular, suppose the initial distribution 
of the process is the Bernoulli distribution $\alpha_0$ 
with marginals 
$$\text{$\alpha_0\{\eta_i=1\}=1-\rho^*$ for $i\le 0$, and  
$\alpha_0\{\eta_i=1\}=\rho^*$ for $i> 0$.}
$$
Let $\alpha_t$ denote the distribution at time $t>0$. Then
any limit point $\mu$ of the time averages $t^{-1}\int_0^t\alpha_s ds$
satisfies {\rm (2.21)} with $\rho=\rho^*$. 
\endproclaim

\hbox{}

\flushpar
{\it 2.3. Remarks and extensions}

\hbox{}

{\it $K$-exclusion and multiple slow bonds.}
Our proof uses the method  of \cite{Se3}. As in that
paper, we could prove Thm.\ 2.2 also
for generalized exclusion processes
that allow $K$ particles per site, instead of just 
one particle. The result would be qualitatively the same,
but with a different growth model in Sect.\ 2.1 and
correspondingly different $\lambda(0)$. 
With suitable large deviation estimates, as in
\cite{Se3}, we envision that the weak law of large 
numbers in Thm.\ 2.2.\ could be strengthened to a 
strong law of large numbers (convergence for almost
every realization of the microscopic evolution). 

A natural extension is to consider finitely many
macroscopically separated slow bonds. Suppose
that for the $n$th process $\eta^n$, jumps
from site $[nx_k]$ occur at rate $\rate_k$
for some set $x_1<\cdots<x_m$ of macroscopic
space points. Thm.\ 2.2  remains unchanged.
The change appears in
 the definition (2.12)
of the macroscopic rate profile: the function $\lambda(\cdot)$ 
 would take on the 
values $\lambda(x_k)=4/\kappa(\rate_k)$, and $\lambda(x)=1$
elsewhere.

\hbox{}

{\it Hamilton-Jacobi equations.}
(2.14) describes $v(x,t)$ as the value function of
an optimal control problem with running cost 
$\lambda(w(s))g_0\bigl({w'(s)}/{\lambda(w(s))}\bigr)$. 
 If 
the function $\lambda$ had sufficient regularity, 
standard theory would imply that $v(x,t)$ is 
 the unique viscosity solution of the
Cauchy problem  
$$v_t + \lambda(x)f_0(v_x)=0\,,\, v(\cdot\,,0)=v_0. 
\tag 2.22
$$
(See Ch.\ 10 in \cite{Ev}.) One can directly check that
the solutions in Cases 1 and 2 of Cor.\ 2.1 satisfy (2.22)
except at jump points. But 
there is currently no theory about existence and 
uniqueness of solutions of Hamilton-Jacobi equations
with a discontinuity of the type that $\lambda(x)$ has
at $0$. The presently available approaches to discontinuous
equations involve smoothing out the discontinuity with
a mollifier, see \cite{Os} and its references. But 
once $\lambda$ has been convolved with a mollifier, the 
discontinuity at $0$ is lost. So a different approach
is required  for our problem.

\hbox{}

{\it The growth model.}
The connection between TASEP and the homogeneous version
of the  growth model of Section 2.1 
 goes back to Rost's 1981 paper
\cite{Ro}. The  last-passage
formulation was first utilized for hydrodynamic limits
and large deviations in \cite{Se1, Se2}. The connection
between exclusion and the growth model 
  will appear  in Section 4. 

An additional motivation for 
the definition of $\kappa(r)$ in terms of the 
path model is this:   analysis of 
models of this type with tools from 
combinatorics and random matrix theory
 has recently led to elegant
exact calculations of limits and fluctuations,
in the work of Baik, Deift, Johansson, and Rains. 
The homogeneous version of this particular growth model is treated
in \cite{Jo}. 

The hydrodynamic result of Thm.\ 2.2 implies a shape
result for the last-passage growth model. We will not
write down the details of such a conversion, but 
only note one point: Suppose the diagonal defect
 does not pass 
through the origin, but instead through the macroscopic
point $(0,u)$ for some  $u> 0$. To capture this, redefine 
the weights in (2.4) by
 $w_{i, i+[nu]}=r^{-1}$, and for other points $w_{i,j}\equiv 1$. 
Now the weights change with $n$, to keep the defect line
at the macroscopic line $y=x+u$. The limiting shape for
this model follows (2.8) part of the way, but develops a kink
around the point where the defect line passes through the 
interface, and the interface is no longer
 convex. This conclusion can be 
worked out from the explicit formulas  of Section 5.1,
via the mapping in Lemma 4.2.

\hbox{}

{\it Relation with the Covert-Rezakhanlou result.}
The definition of 
$\kappa(\rate)$ through 
(2.4)--(2.7) suggests an immediate bound. If we increase
all weights in (2.4) to  $w_{i,j}=1/\rate$ and use
Rost's result (2.8), we get the upper bound 
$\kappa(\rate)\le 4/\rate$. By (2.12) this implies 
$\lambda(0)\ge \rate$, and then from (2.18), 
 $\rho^*\ge (1/2)-(1/2)\sqrt{1-\rate}$. Equivalently, 
the blockage does not disturb a profile at  density $\rho$ if
$\rate\ge 4\rho(1-\rho)$. This is the bound of
Covert and Rezakhanlou \cite{CR}. 

To obtain their bound \cite{CR} approximated the blocked
TASEP with a system whose jump rates depend on spatial
location through a continuous function $\latil(x)$,
so that in the $n$th process jumps from site $i$
occur at rate $\latil(i/n)$. 

The Rost picture shows why a continuous, one-sided {\it macroscopic}
 approximation cannot get any closer than the upper bound 
$\kappa(\rate)\le 4/\rate$. Let $\latil(x,y)$
be a continuous function that satisfies
$\latil(x,y)\le 1$ and $\latil(x,x)=\rate$. 
In (2.5) take new weights $\wtil_{i,j}=\latil(i/n,j/n)^{-1}$, 
and compute $\kappatil(\rate)$ as the limit in (2.7)
with the weights $\wtil_{i,j}$. 
Since $w_{i,j}\le \wtil_{i,j}$, we have $\kappa(\rate)\le
\kappatil(\rate)$. 
Given any $\e>0$, choose $\delta>0$ so that 
$\latil(x,y)<\rate+\e$ in a strip of width $\delta$ 
around the diagonal $\Delta=\{(x,x):0\le x\le 1\}$. 
In the
homogeneous problem (weights $\equiv 1$), the strict concavity 
of the limit (2.8)  implies
that, macroscopically, the diagonal is the unique 
maximizing path for $(x,y)=(1,1)$. Thus with high
probability, the $n\delta$-strip around the microscopic
diagonal $\{(i,i):1\le i\le n\}$ contains a path $\pi$
such that $\sum_\pi Y_{i,j}\ge 4n +o(n)$. The weights
in the $n\delta$-strip  satisfy $\wtil_{i,j}\ge (\rate+\e)^{-1}$,
so it follows
that $\kappatil(\rate)\ge 4(\rate+\e)^{-1}$. Since 
$\e>0$ was arbitrary, $\kappatil(\rate)\ge 4/\rate$,
and  we see that the upper bound
$4/\rate$ is the best we can get by continuously
bounding the rates from below.

\hbox{}

\flushpar
{\bf 3. Proof of Theorem 2.1}

\hbox{}

Let us write $T^\rate_{n,n}$ for
the passage time in (2.5) to indicate dependence on 
$\rate$. Let $0<\rate_1<\rate_2\le 1$. 
Let $\pi$ be the path that gives the maximum in
(2.5) for $T^{\rate_1}_{n,n}$, so that 
$$T^{\rate_1}_{n,n}=\sum_{\mmu\in\pi}Y_\mmu+
\biggl(\frac1{\rate_1}-1\biggr)\sum_{\mmu\in\pi\cap\Delta}Y_\mmu\,.
$$
We wrote $\mmu=(i,j)$ for an integer site on the plane
and $\Delta=\{(x,x):x\in\mmR\}$ is the diagonal. 
Since we are maximizing passage times of paths, 
$$T^{\rate_2}_{n,n}\ge\sum_{\mmu\in\pi}Y_\mmu+
\biggl(\frac1{\rate_2}-1\biggr)\sum_{\mmu\in\pi\cap\Delta}Y_\mmu\,,
$$
and we get
$$T^{\rate_1}_{n,n}-T^{\rate_2}_{n,n}
\le \biggl(\frac1{\rate_1}-\frac1{\rate_2}\biggr)
\sum_{i=1}^n Y_{i,i}\,.
$$
Divide by $n$ and let $n\to\infty$ to get the first inequality
in Thm.\ 2.1. 

Set $\rate_2=1$ and use $\kappa(1)=4$ to get the 
upper bound in (2.9). The value $\kappa(1)=4$ is part
of the lower bound because $\kappa(\cdot)$ is a 
nondecreasing function.   

To get the
other part of the lower bound, consider
this particular path $\pi_n$ from $(1,1)$ to $(n,n)$.
For each $i=1,\ldots, n-1$  choose 
the larger one of $Y_{i+1,i}$ and $Y_{i,i+1}$. 
Since $E[Y'\vee Y'']=3/2$ for two independent rate one
exponentials, this contributes (approximately) $3n/2$ to the sum.
For each pair $\{i,i+1\}$, the above choice led to one
of two situations:

{\it Case 1.} Two points on the same side of the diagonal:
$(i,i+1)$ and  $(i+1,i+2)$, or $(i+1,i)$ and  $(i+2,i+1)$.

{\it Case 2.}
Two points on opposite sides of the diagonal:
$(i,i+1)$ and  $(i+2,i+1)$, or $(i+1,i)$ and  $(i+1,i+2)$.

To complete the path 
$\pi$  we pick one site between each pair chosen above. 
  In {\it Case 1}
we can choose the larger one of the diagonal $r^{-1}Y_{i+1,i+1}$
and an off-diagonal value, $Y_{i,i+2}$ or  $Y_{i+2,i}$
depending on the subcase of {\it Case 1}. The expected
contribution is 
$E[Y'\vee (r^{-1}Y'')]= (\rate^2+\rate+1)\rate^{-1}(1+\rate)^{-1}$. 
 In {\it Case 2} there
is no choice: the path must go through $(i+1,i+1)$ and 
take the diagonal value $r^{-1}Y_{i+1,i+1}$. 
 {\it Cases 1} and {\it 2}  are equally likely, so this 
second step contributes approximately
$(n/2)\cdot  (\rate^2+\rate+1)\rate^{-1}(1+\rate)^{-1}
 + (n/2)\cdot \rate^{-1}$. Adding the contributions 
from the two steps gives the lower bound in (2.9).

\hbox{}

\flushpar
{\bf 4. Proof of Theorem 2.2}

\hbox{}

\flushpar
{\it 4.1. Construction of the process and the variational 
coupling}

\hbox{}

We follow the construction in Section 4 of \cite{Se3},
and  only outline it here. 
We construct a process $z(t)= (z_i(t):i\in\mmZ)$ 
 of labeled particles that move on  $\mmZ$. 
The location of the $i$th particle at time $t$ is 
$z_i(t)$, and these satisfy 
$$0\le z_{i+1}(t)-z_i(t)\le 1\,.
\tag 4.1
$$
For the dynamics, let $\{\DD_i\}$ be a collection of mutually
independent Poisson jump time processes on the time line 
$(0,\infty)$. $\DD_0$ has rate $\rate$, and all other 
$\DD_i$ have rate 1. In the graphical construction, 
 $z_i$ attempts a jump one step to the {\it left} at 
epochs of $\DD_i$. The jump is executed if it does not
violate (4.1). This happens independently for all $i$. 

Once the $z(t)$ process is constructed, the exclusion
process $\eta(t)$ is defined by
$$\eta_i(t)=z_i(t)-z_{i-1}(t)\,.
\tag 4.2
$$ 
The $z$-particles keep track of the current of the 
exclusion process $\eta(\cdot)$. 
The number of $\eta$-particles that have 
left site $i$ in time $[0,t]$ is given by 
$$J_i(t)=z_i(0)-z_i(t)\,.
\tag 4.3
$$

Assume now that the  process 
$z(\cdot)$ has been constructed on some probability space
that supports  the initial
configuration $z(0)=(z_i(0))$, and the Poisson processes
$\{\DD_i\}$ that are  independent of $(z_i(0))$. 
We define a family $\{w^k:k\in\mmZ\}$ of auxiliary processes 
on this same probability space. Each  
$w^k(t)=(w^k_i(t):i\in\mmZ)$ is an exclusion process
just like  $z(t)$, so  (4.1) is in force for all $i$. 
 Initially 
$$w^k_i(0)=\cases
z_k(0)\,, &i\ge 0\\
z_k(0)+i\,, &i<0\,
\endcases
\tag 4.4
$$
and dynamically
$$\text{$w^k_i$ attempts to jump to $w^k_i-1$ at the epochs of $\DD_{i+k}$.}
\tag 4.5
$$
The usefulness of the family
of processes $\{w^k\}$ lies in this fact:

\proclaim{Lemma 4.1} For all $i\in\mmZ$ and $t\ge 0$,
$$z_i(t)=\sup_{k\in\mmZ} w^k_{i-k}(t)\ \ \text{a.s.}
\tag 4.6
$$
\endproclaim

This lemma is proved as Lemma 4.2 in \cite{Se3} so we will not
repeat the proof here. This is the ``variational coupling''
that is the key to our proof. 

Thinking of the process $w^k(t)$ as giving the height
of an  interface over the sites  $i$, we want to normalize
it to start at height zero, and also to advance in the 
increasing coordinate 
direction. Hence define a new family of processes
$\{\xi^k\}$ by 
$$\xi^k_i(t)=z_k(0)-w^k_i(t)\ \ \text{ for $i\in\mmZ$, $t\ge 0$.}
\tag 4.7
$$
Now we can write (4.6) as 
$$
z_i(t)=\sup_{k\in\mmZ} \{z_k(0)-\xi^k_{i-k}(t)\}\,.
\tag 4.8
$$
The virtue of (4.8) is that the effect of the initial
condition $(z_k(0))$ has been separated from the effect
of the Poisson jump times $\{\DD_i\}$. The process 
$\xi^k$ does not depend on $z_k(0)$, and depends on 
the superscript $k$ only through a translation of the
indexing of 
$\{\DD_i\}$. 
 Initially 
$$\xi^k_i(0)=\cases
0\,, &i\ge 0\\
-i\,, &i<0\,.
\endcases
\tag 4.9
$$
 The dynamical rule for the $\xi^k$ process
is that 
$$
\text{$\xi^k_i$ jumps to $\xi^k_i+1$ at epochs of $\DD_{i+k}$,}
\tag 4.10
$$
provided the inequalities 
$$\xi^k_{i}\le \xi^k_{i-1}\qquad\text{and}\qquad
\xi^k_{i}\le \xi^k_{i+1}+1
\tag 4.11
$$
are not violated. Notice that $\xi^k_{-k}$ jumps at rate $\rate$,
while other $\xi^k_i$ jump at rate 1. 

We can now outline the strategy for proving Thm.\ 2.2.
Given the initial configurations $\eta^n(0)=(\eta^n_i(0):i\in\mmZ)$
that appear in hypothesis (2.16), define initial
configurations $z^n(0)=(z^n_i(0):i\in\mmZ)$ so that 
$z^n_0(0)=0$ and (4.2) holds at time $t=0$. 
Then hypothesis (2.16) implies that 
$$\lim_{n\to\infty} n^{-1}z^n_{[nq]}(0)=v_0(q)
\ \ \text{in probability}
\tag 4.12
$$
for all $q\in\mmR$, with $v_0$ defined by (2.11). 

Construct the
processes $z^n(t)$ as indicated above, and define 
the exclusion processes $\eta^n(t)$ by (4.2). 
Define $v(x,t)$ by (2.14). 
From (4.2)--(4.3) we see that both limits of Thm.\ 2.2
follow from proving that for all
$x\in\mmR$ and $t>0$, 
$$\lim_{n\to\infty} n^{-1}z^n_{[nx]}(nt)=v(x,t)
\ \ \text{in probability.}
\tag 4.13
$$

Now rewrite (4.8) with the correct scaling:
$$
n^{-1}z^n_{[nx]}(nt)=\sup_{q\in\mmR} \lbrakk n^{-1}z^n_{[nq]}(0)-
n^{-1}\xi^{[nq]}_{[nx]-[nq]}(nt)\rbrakk\,.
\tag 4.14
$$
In this formula each process $z^n(\cdot)$ is defined on
a probability space that supports the initial 
configuration $z^n(0)$ and the Poisson processes
$\{\DD_i\}$. On each such probability space we define
the processes $\{\xi^k(\cdot)\}$ as functions of $\{\DD_i\}$,
according to (4.9)--(4.11). 
The proof of (4.13) is now to show that the right-hand 
side of (4.14) converges to the right-hand 
side of (2.14). The first step is to prove a limit 
for the $\xi$-term. 

\hbox{}

\flushpar
{\it 4.2. Limit for $\xi$}

\hbox{}

First some definitions. The meaning of these notions 
will be explained below.  Let
$$\VV=\{(x,y)\in\mmR^2: y\ge 0\,,\, x\ge -y\}\,.
$$
For $(x,y)\in\VV$, let
$$\gamma_0(x,y)=\bigl(\sqrt{x+y\,}+\sqrt{y\,}\bigr)^2\,.
\tag 4.15
$$
Let $\mmx(s)=(x_1(s),x_2(s))$ denote a path in $\mmR^2$,
defined on some interval of $s$-values. 
 For $(x,y)\in\VV$ and $q\in\mmR$ let
$$\Gamma^q(x,y)=\sup\lbrakk \int_0^1 \frac
{\gamma_0\bigl(\mmx'(s)\bigr)}{\lambda\bigl(x_1(s)-q\bigr)}\,ds: 
\mmx(\cdot)\in\XX(x,y)\rbrakk\,,
\tag 4.16
$$
where $\XX(x,y)$ is the collection of piecewise 
$C^1$ paths $\mmx:[0,1]\to\VV$ that satisfy 
$$\text{$\mmx(0)=(0,0)$, $\mmx(1)=(x,y)$, 
 and $\mmx'(s)\in\VV$  for all $s$.}
\tag 4.17
$$
The last condition ensures that 
 $\gamma_0\bigl(\mmx'(s)\bigr)$ is defined. 
The function $\lambda(\cdot)$ in the definition of $\Gamma^q$
is the macroscopic rate defined by (2.12). 
Lastly, for  $q,x\in\mmR$ and $t>0$, set 
$$g^q(x,t)=\inf\{ y: (x,y)\in\VV, \Gamma^q(x,y)\ge t\}\,.
\tag 4.18
$$
$\Gamma^q(x,y)$ represents the macroscopic time it takes
a $\xi$-type interface process to reach point $(x,y)$. The
point $q$ marks the $x$-coordinate of the defect column that 
(macroscopically) has the slow rate $\lambda(0)$. 
The level curve of $\Gamma^q$ given by $g^q(\cdot\,,t)$ 
represents the limiting interface of a $\xi$-process, as 
stated in the next proposition. 

\proclaim{Proposition 4.1} For all $q,x\in\mmR$ and $t>0$, 
$$\lim_{n\to\infty} n^{-1}\xi^{[nq]}_{[nx]}(nt)=g^{-q}(x,t)
\ \ \text{in probability.}
\tag 4.19
$$
\endproclaim

\demo{Remark 4.1} For the reader familiar with the proof of
\cite{Se3}, let us point out that Prop.\ 4.1 is the 
analogue of Cor.\ 5.1 in \cite{Se3}, and 
$\Gamma^q$ corresponds to the limit in 
Prop.\ 5.1 in \cite{Se3}. 
\enddemo

To prove Prop.\ 4.1 we follow Section 5
of \cite{Se3} and switch to a last-passage representation.
Lemma 4.2 below shows  how the growth model of Section 2.1
enters the picture, and justifies (2.7) as the definition
of $\kappa(\rate)$. 
 
Define a lattice analogue of the interior of the wedge $\VV$ by
$\LL=\{(i,j)\in\mmZ^2: j\ge 1, i\ge -j+1\}$, with
boundary
$\partial\LL=\{(i,0): i\ge 0\} \cup\{(i,-i): i<0\}$. 
For $(i,j)\in\LL\cup\partial\LL$, let
$$L^k(i,j)=\inf\{t\ge 0: \xi^k_i(t)\ge j\}
\tag 4.20
$$
denote the time when $\xi^k_i$ reaches level $j$. 
The rules (4.9)--(4.11) give the boundary conditions
$$\text{$L^k(i,j)=0$ for $(i,j)\in\partial\LL$,}
\tag 4.21
$$
and for $(i,j)\in\LL$ the equation 
$$L^k(i,j)=\max\{ L^k(i-1,j), L^k(i,j-1), L^k(i+1,j-1)\}
+\beta^k_{i,j}\,,
\tag 4.22
$$
where $\beta^k_{i,j}$ is an exponential waiting time, 
independent of everything else. It represents the time
$\xi^k_i$ waits to jump, {\it after} $\xi^k_i$ and its
neighbors $\xi^k_{i-1}$, $\xi^k_{i+1}$ have reached 
the positions that permit
$\xi^k_i$ to jump to $j$. By (4.10),  $\beta^k_{-k,j}$
has rate $\rate$, but for $i\ne -k$, $\beta^k_{i,j}$ has rate 1. 

The waiting time  $\beta^k_{i,j}$ cannot be read directly 
from $\DD_{i+k}$. One has to construct the evolution
$\xi^k(\cdot)$ 
up to the stopping time 
$L=\max\{ L^k(i-1,j), L^k(i,j-1), L^k(i+1,j-1)\}$, and then 
$\beta^k_{i,j}$ is the 
waiting time to the next epoch in $\DD_{i+k}$ after $L$. 
The last-passage representation entails switching
probability spaces so that the waiting times 
$\beta^k_{i,j}$ become the basic building blocks of the 
construction. We shall switch notation to keep the 
two constructions distinct. 

We now construct a last-passage
 growth model on $\LL$ that 
has a  defect in the column $\{(m,j): j\in\mmZ\}$, where $m\in\mmZ$
is fixed. 
Let 
$\{\tau_{\mmu}:\mmu\in\LL\}$ denote a collection of i.i.d.\ 
exponential rate 1 random variables. Define weights 
$$\omega^{m}_{i,j} =\cases 1\,, &i\ne m\\
1/\rate\,, &i= m\,.
\endcases
\tag 4.23
$$
Given $\mmu=(u,v)\in\LL$, let $\Pi(\mmu)$ denote the set of
lattice paths 
$\pi=\{(0,1)=(i_1,j_1),(i_2,j_2), \ldots, (i_p,j_p)=\mmu\}$
whose admissible steps satisfy
$$\text{$(i_\ell,j_\ell)-(i_{\ell-1},j_{\ell-1})=(1,0)$, 
$(0,1)$, or $(-1,1)$.}
\tag 4.24
$$
Finally, define the passage times $\{T^{m}(\mmu)\}$
by
$$\text{$T^{m}(\mmu)=0$ for $\mmu\in\partial\LL$,}
\tag 4.25
$$
and for $\mmu\in\LL$ as the maximal weighted sum of waiting times
over admissible paths:
$$T^{m}(\mmu)=\max_{\pi\in\Pi(\mmu)}
\sum_{\mmv\in\pi}\omega^{m}_\mmv\tau_\mmv\,.
\tag 4.26
$$
Comparison of (4.21)--(4.22) and (4.25)--(4.26) shows that the passage 
time processes $\{L^k(\mmu):\mmu\in\LL\cup\partial\LL\}$ and 
$\{T^{-k}(\mmu):\mmu\in\LL\cup\partial\LL\}$ have the same 
distribution. The superscript changes 
by a minus sign from $L^k(\mmu)$ to $T^{-k}(\mmu)$
because $L^k(\mmu)$ is the passage time for the process
$\xi^k$, but the slow column for this process is at $-k$, so 
we set $m=-k$ in  $T^{m}(\mmu)$. 

The conclusion of this is that $T^{-k}(i,j)$ has the 
same distribution as the time when $\xi^k_i$ reaches level $j$. 
So, as in \cite{Se3}, Prop.\ 4.1 will follow if we prove 

\proclaim{Proposition 4.2} For all $q\in\mmR$ and $(x,y)$
in the interior of $\VV$, 
$$\lim_{n\to\infty} n^{-1}T^{[nq]}([nx],[ny])=\Gamma^{q}(x,y)
\ \ \text{in probability.}
\tag 4.27
$$
\endproclaim

\demo{Remark 4.2}
Let us first explain the known homogeneous situation
where $\rate=1$ and $\lambda(0)=1$.
Write $\xi^{(\rate=1)}$ and $T^{(\rate=1)}$ for the
interface process defined by (4.9)--(4.11)
and the passage time in (4.26) when  $\rate=1$
and there 
is no special column. 
 In this case the 
limits in (4.19) and  (4.27) are  given by
$$\aligned
&\lim_{n\to\infty} n^{-1}\xi^{(\rate=1)}_{[nx]}(nt)=tg_0(x/t)
= \frac{t}4\biggl(1-\frac{x}t\biggr)^2\quad
\text{ for $-t\le x\le t$, and }\\
&\lim_{n\to\infty} n^{-1}T^{(\rate=1)}([nx],[ny])=\gamma_0(x,y)
= \bigl(\sqrt{x+y\,}+\sqrt{y\,}\bigr)^2\,.
\endaligned
\tag 4.28
$$
 The explicit values cannot be 
inferred directly from the path model, but 
indirectly via the connection with TASEP. Briefly, 
one computes the current $f_0(\rho)=\rho(1-\rho)$  
from the known invariant distributions of TASEP.
The coupling (4.8) implies that the limiting shape $g_0$
for $\xi$ is the conjugate of $f_0$, and one can derive
 $g_0$. The 
shape $g_0$ is a level curve of the passage time
$\gamma_0$, so one obtains the formula for $\gamma_0$ 
from 
$$\gamma_0(x,g_0(x))=1\,,
\tag 4.29
$$ 
and the homogeneity
of $\gamma_0$. This type of argument is repeated
in the examples in \cite{Se1}. 
\enddemo

\demo{Proof of Proposition 4.2} Consider the 
possible maximizing macroscopic curves in (4.16).  By the concavity of
$\gamma_0$, of all the paths
$\mmx(s)=(x_1(s),x_2(s))$ in $\XX(x,y)$, we only
need to consider these two types:
$$\mmx(s)=(sx,sy)\,,\ \ 0\le s\le 1
\tag 4.30
$$
with value [except in the case $x=q=0$ which is covered by (4.32)]
$$\int_0^1 \frac
{\gamma_0\bigl(\mmx'(s)\bigr)}{\lambda\bigl(x_1(s)-q\bigr)}\,ds
=\gamma_0(x,y)\,,
\tag 4.31
$$
and 
$$\mmx(s)=\cases \frac{s}{s_1}(q, x_2(s_1))\,, &0\le s\le s_1\\
(q,x_2(s_1))+\frac{s-s_1}{s_2-s_1}(0, x_2(s_2)-x_2(s_1))\,, 
&s_1\le s\le s_2\\
(q,x_2(s_2))+\frac{s-s_2}{1-s_2}(x-q, y-x_2(s_2))\,, 
&s_2\le s\le 1\,,
\endcases
\tag 4.32
$$
with value
$$\aligned
&\int_0^1 \frac
{\gamma_0\bigl(\mmx'(s)\bigr)}{\lambda\bigl(x_1(s)-q\bigr)}\,ds\\
=\;&\gamma_0\bigl(\mmx(s_1)\bigr)+\lambda(0)^{-1}
\gamma_0\bigl(\mmx(s_2)-\mmx(s_1)\bigr)+
\gamma_0\bigl((x,y)-\mmx(s_2)\bigr)\,.
\endaligned
\tag 4.33
$$
In words: the path of type (4.30) is a single line segment
from $(0,0)$ to $(x,y)$. The path of type (4.32) first uses
parameter interval $[0,s_1]$ to take 
a straight line path to the vertical line $x=q$, then
spends  interval $[s_1,s_2]$ on this line
 to take advantage
of the slow rate $\lambda(0)$, 
and finally takes a straight line path
to $(x,y)$. 

Now consider the microscopic path problem (4.26). We need to
establish the connection between it and the quantity
$\lambda(0)=4/\kappa(\rate)$ defined in Section 2.1. 

\proclaim{Lemma 4.2} Set $q=0$ so that the special 
column goes through the origin. Then for $y>0$
$$\lim_{n\to\infty} n^{-1}T^{0}(0,[ny])=\kappa(\rate)y
\ \ \text{in probability.}
\tag 4.34
$$
\endproclaim

\demo{Proof of Lemma 4.2} The Lemma follows because 
the growth model of Sect.\ 2.1 with a diagonal defect is the same 
as the one studied here when the columnar defect is at the 
origin. A simple mapping 
reveals this. First observe that 
the admissible step $(0,1)$ can be eliminated from (4.24), 
 because each $(0,1)$-step in a path
 can be  replaced by a $(1,0)$-step
followed by a $(-1,1)$-step. This change adds 
a site to the path and hence increases its overall
passage time. So we may assume that
 $\Pi(\mmu)$ contains only paths that have 
 admissible steps $(1,0)$ and $(-1,1)$. 

Consider the bijection $\psi:\LL\to\mmN^2$ given by 
$\psi(i,j)=(i+j,j)$. For $(i,j)\in\mmN^2$, 
define random variables $Y_{i,j}=
\tau_{\psi^{-1}(i,j)}$ and weights $w_{i,j}=
\omega^0_{\psi^{-1}(i,j)}$. Then $w_{i,j}$ satisfies (2.4). 
For $\pi\in\Pi(0,n)$, the image path $\psi(\pi)$ 
runs from $(1,1)$ to $(n,n)$, and has steps of
two kind: $(1,0)$ and $(0,1)$. Thus the map $\psi$
transforms  $T^0(0,n)$ of (4.26) into $T_{n,n}$ of (2.5). 
Now the lemma follows from the definition (2.7) of
$\kappa(\rate)$. \qed
\enddemo

Return to the proof of Prop.\ 4.2. To first prove 
$$\liminf_{n\to\infty} n^{-1}T^{[nq]}([nx],[ny])\ge\Gamma^{q}(x,y)\,,
\tag 4.35
$$
consider any macroscopic
path $\mmx(\cdot)$ of type (4.32). [We leave the 
easier type (4.30) to the reader.] Let $\pi_n$ be the 
microscopic path through the sites $(0,1)$, 
$[n\mmx(s_1)]\equiv([nx_1(s_1)],[nx_2(s_1)])$,
$[n\mmx(s_2)]$, and $([nx],[ny])$, constructed 
so that each of the three segments maximizes 
passage time between its endpoints. Then by (4.34) and the limit
(4.28) for the 
homogeneous case, 
$$\aligned
&\liminf_{n\to\infty} n^{-1}T^{[nq]}([nx],[ny])\ge
\liminf_{n\to\infty}n^{-1}\sum_{\mmv\in\pi_n}\omega^{[nq]}_\mmv
\tau_\mmv \\
\ge\;&\gamma_0\bigl(\mmx(s_1)\bigr)+
\kappa(\rate)\bigl(x_2(s_2)-x_2(s_1)\bigr)+ 
\gamma_0\bigl((x,y)-\mmx(s_2)\bigr)\\
=\;&\text{the value of the path in (4.33).}
\endaligned
\tag 4.36
$$
The last equality follows because 
$\lambda(0)^{-1}
\gamma_0\bigl(\mmx(s_2)-\mmx(s_1)\bigr)
=(1/4)\kappa(\rate)\gamma_0\bigl(0,x_2(s_2)-x_2(s_1)\bigr) 
=\kappa(\rate)\bigl(x_2(s_2)-x_2(s_1)\bigr)$. 
The reason there might not be equality in the last
inequality in (4.36) is that an optimal path between, say, 
$(0,1)$ and $[n\mmx(s_1)]$ might actually take advantage
of the $[nq]$-column and return a larger value
than $\gamma_0\bigl(\mmx(s_1)\bigr)$. 
A similar argument for paths of type (4.30) justifies (4.35). 

Now for  the complementary upper bound 
$$\limsup_{n\to\infty} n^{-1}T^{[nq]}([nx],[ny])\le\Gamma^{q}(x,y)\,.
\tag 4.37
$$
 Each macroscopic path 
$\mmx(\cdot)$ in $\XX(x,y)$ is contained in a fixed compact
subset $A$ of $\VV$.  
Choose  $\delta>0$ so that   
$$\text{$|\gamma_0(\mmx)-\gamma_0(\mmy)|< \e $
for $\mmx,\mmy\in A$ such that $|\mmx-\mmy|<\delta$,}
\tag 4.38
$$
and then a partition
$$0=b_0<b_1<\cdots<b_s=y$$
of $[0,y]$ 
with mesh $\max (b_{i+1}-b_i) <\delta$. Let $\mmx^{(i,j)}$
be the path of type (4.32) with $x^{(i,j)}_2(s_1)=b_i<
b_j=x^{(i,j)}_2(s_2)$. 

Let us adopt the following 
 generalization of the notation in (4.26):
$T^m(\mmu,\mmv)$ denotes the maximal weighted sum over
admissible paths from $\mmu$ to $\mmv$. So
$T^m(\mmu)$ in (4.26) is the same as $T^m((0,1),\mmu)$. 

Let $\pi_n$ be the maximizing microscopic path in (4.26) for
$\mmu=([nx],[ny])$ and $m=[nq]$. The easy situation is when 
$\pi_n$ does not intersect the vertical column 
$([nq],j)$ that has the slow rate $r$. Then $T^{[nq]}([nx],[ny])$
equals the homogeneous passage time $T^{(\rate=1)}([nx],[ny])$.
If this happens infinitely often 
along the subsequence taken on the left-hand side
of (4.37), then (4.37) follows from (4.28). 

Otherwise, pick indices $k\le \ell$
such that the path $\pi_n$ first touches the vertical column 
$([nq],j)$ in the range $[nb_k]\le j< [nb_{k+1}]$, and 
for the last time in the range $[nb_\ell]\le j< [nb_{\ell+1}]$. 
Then quite obviously 
$$\aligned
&T^{[nq]}([nx],[ny])=\sum_{\mmv\in\pi_n}\omega^{[nq]}_\mmv
\tau_\mmv\\
\le\; &T^{(\rate=1)}\bigl([nq],[nb_{k+1}]\bigr) 
+T^{[nq]}\bigl(([nq],[nb_{k}]),([nq],[nb_{\ell+1}])\bigr)\\
&\qquad\qquad 
+T^{(\rate=1)}\bigl(([nq],[nb_{\ell}]),([nx],[ny])\bigr)\\
\le\; &\max_{i\le j}\lbrakk 
T^{(\rate=1)}\bigl([nq],[nb_{i+1}]\bigr) 
+T^{[nq]}\bigl(([nq],[nb_{i}]),([nq],[nb_{j+1}])\bigr)\\
&\qquad\qquad 
+T^{(\rate=1)}\bigl(([nq],[nb_{j}]),([nx],[ny])\bigr)
\rbrakk\,.
\endaligned
$$
Divide by $n$, let $n\to\infty$,  use the limits 
 (4.28) and (4.34), and then (4.38) 
to get, with $C=1+\lambda(0)^{-1}$, 
$$\aligned
&\limsup_{n\to\infty} n^{-1}T^{[nq]}([nx],[ny])\\
\le\;&\max_{i\le j}\lbrakk
\gamma_0(q,b_{i+1})+\kappa(\rate)(b_{j+1}-b_i)+
\gamma_0(x-q,y-b_j)\rbrakk\\
\le\;&\max_{i\le j}\lbrakk
\gamma_0(q,b_{i})+\lambda(0)^{-1}\gamma_0(0, b_j-b_i) +
\gamma_0(x-q,y-b_j)\rbrakk+C\e\\
=\;&\max_{i\le j}\lbrak
 \text{value of path $\mmx^{(i,j)}(\cdot)$} \rbrak+C\e\\
\le\;&\Gamma^q(x,y)+C\e\,.
\endaligned
$$
This completes the proof of Proposition 4.2.
\qed
\enddemo

Prop.\ 4.1 follows from Prop.\ 4.2 as Cor.\ 5.1 follows
from Prop.\ 5.1 in \cite{Se3}. 

\hbox{}

\flushpar
{\it 4.3. Hydrodynamic limit}

\hbox{}

Using (4.12), (4.14), and (4.19) we can now prove that,
in probability, 
$$\lim_{n\to\infty} n^{-1}z^n_{[nx]}(nt)=
\vtil(x,t)\equiv\sup_{q\in\mmR}\lbrak v_0(q)-g^{-q}(x-q,t)\rbrak\,.
\tag 4.39
$$
The argument is the one  from eqn.\ (6.4) to (6.15) in \cite{Se3}, 
so we will not repeat it here. To complete the proof of
(4.13) and thereby the proof of Thm.\ 2.2, we need to
show that  the limiting value $\vtil(x,t)$
defined above  agrees with the desired limit
$v(x,t)$ defined by (2.14). 

\hbox{}

\flushpar
{\it 4.4. Formula for $v(x,t)$}

\hbox{}

By (4.16)--(4.18), the definition (4.39) 
of $\vtil(x,t)$ can be rewritten as
$$\aligned
\vtil(x,t)
&=\sup_{q,y\in\mmR}\lbrakk v_0(q)-y: \text{there exists
a path $\mmx(\cdot)\in\XX(x-q,y)$}\\
&\qquad\qquad\qquad \text{such that }\ \ 
\int_0^1 \frac
{\gamma_0\bigl(\mmx'(s)\bigr)}{\lambda\bigl(x_1(s)+q\bigr)}\,ds\ge t
\rbrakk\,.
\endaligned
\tag 4.40
$$

\proclaim{Proposition 4.3} $v(x,t)=\vtil(x,t)$. 
\endproclaim

\demo{Proof} The proof involves mapping the paths 
$\mmx(\cdot)$  
to the paths $w(\cdot)$ that appear in (2.14), and vice versa. 

Given a path $\mmx(\cdot)\in\XX(x-q,y)$
that appears in (4.40), define a new time
variable $\tau=\tau(s)$ by 
$$\tau(s)=\int_0^s \frac
{\gamma_0\bigl(\mmx'(s)\bigr)}{\lambda\bigl(x_1(s)+q\bigr)}\,ds\,.
\tag 4.41
$$
Let the terminal $\tau$-time be  $t_1=\tau(1)$. From 
(4.40) we know that $t_1\ge t$. Let $s=s(\tau)$ be the 
inverse time change. Define a path $z:[0,t_1]\to\mmR$ by 
$$z(\tau)=q+x_1(s(\tau))\,.
\tag 4.42
$$
Then 
$$z(0)=q\,,\, z(t_1)=x\,,\,\text{ and }\ z'(\tau)=x_1'(s(\tau))s'(\tau)\,.
\tag 4.43
$$
Differentiating (4.41), relation (4.29), and the homogeneity
of $\gamma_0$ [means: $\gamma_0(cx,cy)=c\gamma_0(x,y)$] give
$$x_2'(s)=\tau'(s)\lambda\bigl(x_1(s)+q\bigr)
g_0\biggl(\frac{x_1'(s)}{\tau'(s)\lambda\bigl(x_1(s)+q\bigr)}\biggr)\,.
\tag 4.44
$$
Since $y=x_2(1)$, we can use this to compute
$$\aligned
y&=\int_0^1 x_2'(s)ds=\int_0^1 \tau'(s)\lambda\bigl(x_1(s)+q\bigr)
g_0\biggl(\frac{x_1'(s)}{\tau'(s)\lambda\bigl(x_1(s)+q\bigr)}\biggr)ds\\
&=\int_0^{t_1} \lambda(z(\tau))g_0\biggl(\frac{z'(\tau)}{\lambda(z(\tau))}\biggr)
d\tau\,.
\endaligned
$$
The only problem is that $z(\cdot)$ is defined on $[0,t_1]$
instead of on the possibly smaller interval $[0, t]$. 
Let $w(\sigma)=z(t_1\sigma/t)$ be a time change of $z(\cdot)$
defined for $0\le\sigma\le t$. Then, because 
$sg_0(x/s)$ is nondecreasing in $s$, change of variable 
$\tau=(t_1/t)\sigma$ in the last integral above shows 
that
$$v_0(q)-y\le v_0(w(0))- \int_0^{t} \lambda(w(\sigma))g_0\biggl(\frac{w'(\sigma)}
{\lambda(w(\sigma))}\biggr)
d\sigma\,.
$$
Since $\mmx(\cdot)$ was an arbitrary path inside the braces in (4.40), 
we have shown that $\vtil(x,t)\le v(x,t)$. 

Conversely, take a path $w:[0,t]\to\mmR$ that appears
in (2.14). 
Define  a path $\mmx(\sigma)=(x_1(\sigma),x_2(\sigma))$ 
for $\sigma\in[0,1]$ by
$$\aligned
x_1(\sigma)&=w(\sigma t)-w(0)\\
x_2(\sigma)&=\int_0^{\sigma t} \lambda(w(s))g_0\biggl(\frac{w'(s)}
{\lambda(w(s))}\biggr)
ds\,.
\endaligned
\tag 4.45
$$
 Let $q=w(0)$ and $y=x_2(1)$. Then 
$\mmx(\cdot)\in\XX(x-q,y)$, provided
$\mmx'(\sigma)\in\VV$ [recall conditions (4.17)].
This follows because
$\lambda$ and $g_0$ are nonnegative functions, and 
because $g(z)\ge -z$ for all $z$. 
Also, (4.29) and (4.45) give 
$$ \frac
{\gamma_0\bigl(\mmx'(\sigma)\bigr)}{\lambda\bigl(x_1(\sigma)+q\bigr)}= t\,,
$$
so the integral condition inside the braces in (4.40) is
satisfied. We conclude that $\mmx(\cdot)$ is a
 path that appears in (4.40), and since
$$v_0(w(0))-\int_0^{t} \lambda(w(s))g_0\biggl(\frac{w'(s)}
{\lambda(w(s))}\biggr)
ds=v_0(q)-y\le \vtil(x,t)\,,$$
we have 
$v(x,t)\le \vtil(x,t)$.  This completes the proof of 
Prop.\ 4.3.
\qed
\enddemo

We have now proved Thm.\ 2.2.

\hbox{}

\flushpar
{\bf 5. Proofs of the Corollaries}

\hbox{}

\flushpar
{\it 5.1. Proof of Corollary 2.1 }

\hbox{}

As a preliminary step for calculating macroscopic profiles
from the variational 
formula (2.14), we optimize the integral term as a 
function of the initial point  $q=w(0)$. So let 
$$\aligned
I(x,t,q)=&\inf\lbrakk \int_0^t \lambda(w(s))
g_0\biggl(\frac{w'(s)}{\lambda(w(s))}\biggr)
ds: \\
&\text{$w:[0,t]\to\mmR$ is 
piecewise $C^1$, $w(0)=q$, and  $w(t)=x$} \rbrakk\,.
\endaligned
\tag 5.1
$$
In terms of the limiting shapes of Prop.\ 4.1, 
$I(x,t,q)=g^{-q}(x-q,t)$, so from the formulas below the 
reader can deduce explicit expressions for the limits in 
(4.19). 

By the convexity of $g_0$, it suffices to consider 
the following two types of paths in (5.1): 
either $w(\cdot)$ is a single linear segment from 
$q= w(0)$ to $w(t)=x$; or it consists of a linear
segment from 
$q= w(0)$ to $w(s_1)=0$, a constant segment
$w(s)=0$ for $s_1\le s\le s_2$, and a linear segment
from $w(s_2)=0$ to  $w(t)=x$.
Only calculus is involved in finding the optimal
paths, so we skip the details and present a summary
of the results. 
Abbreviate
$$B=\sqrt{1-\lambda(0)\,}\,.$$

Five different ranges of the 
variables $x$ and $q$ appear.
$$|x|\ge Bt \,\,; 
\tag 5.2a
$$
$$
\aligned
&\text{$0\le x< Bt $ , and }\\
&\text{$q\le x-Bt $ or 
$q\ge \bigl( \sqrt{Bt\,}-\sqrt{x}\,\bigr)^2$}\,\,;
\endaligned
\tag 5.2b
$$
$$\aligned
&\text{$ -Bt <x\le 0$ , and }\\
&\text{$q\le -\bigl( \sqrt{Bt\,}-\sqrt{|x|}\,\bigr)^2$ or 
$q\ge x+Bt $}\,\,;
\endaligned
\tag 5.2c
$$
$$
\aligned
&\text{$0\le x< Bt $ , and }\\
&\text{$ x-Bt <q
<\bigl( \sqrt{Bt\,}-\sqrt{x}\,\bigr)^2$}\,\,;
\endaligned
\tag 5.2d
$$
$$\aligned
&\text{$ -Bt <x\le 0$ , and }\\
&\text{$-\bigl( \sqrt{Bt\,}-\sqrt{|x|}\,\bigr)^2<q
< x+Bt $}\,\,.
\endaligned
\tag 5.2e
$$

These are the optimal values: 

\hbox{}

Cases (5.2a--c): $I(x,t,q)=tg_0\bigl((x-q)/t\bigr)$ and 
the optimal path is $w(s)=q+(s/t)(x-q)$.

\hbox{}

Cases (5.2d--e): 
$$\aligned
I(x,t,q)=&\frac{|q|}{B} 
g_0\bigl(-B{q}/{|q|}\,\bigr)
+\biggl( t-\frac{|x|+|q|}{B}\biggr) \lambda(0)g_0(0)
+\frac{|x|}{B}
 g_0\bigl(B{x}/{|x|}\,\bigr)\\
=&\frac{B}2\bigl({|x|+|q|}\bigr)-\frac{x-q}2+\frac{t}4\lambda(0)
\endaligned
$$
and 
the optimal path is 
$$w(s)=\cases
q-sq/s_1\,, &0\le s<s_1\\
0\,, &s_1\le s<s_2\\
(s-s_2)x/(1-s_2)\,, &s_2\le s\le 1\,,
\endcases
$$
with 
$$s_1={|q|}/{B}\qquad\text{and}\qquad
s_2=t-{|x|}/{B}\,.
$$

\hbox{}

Proof of Corollary 2.1 is now reduced to finding
$$v(x,t)=\sup_{q\in\mmR}\{v_0(q)-I(x,t,q)\}
\tag 5.3
$$
with $v_0(q)=\rho q$. We skip the calculus details.

\hbox{}

\flushpar
{\it 5.2. Proof of Corollary 2.2 }

\hbox{}

Assume now that $\rho^*<1/2$, in other words, that the 
slow bond disturbs the hydrodynamic profiles. 
First check from (5.3) that all the profiles admitted 
by the restrictions stated in Cor.\ 2.2 are in fact 
invariant. This contains the following cases: constants
at $\rho^*$ and at $1-\rho^*$; a piecewise constant
 profile with 
a single entropy shock [jump from $\rho^*$ to $1-\rho^*$]
anywhere in $\mmR$; a piecewise constant
 profile with 
a single non-entropy shock [jump from $1-\rho^*$ to $\rho^*$]
at $x=0$;  and a piecewise constant
 profile with 
a non-entropy shock at $x=0$, and  an entropy shock
in $(-\infty,0)$, or in $(0,\infty)$, or in both. 

To prove that these are the only possible invariant profiles, 
let $\rho_0$ be a profile such that $\rho^*\le\rho_0(x)\le 1-\rho^*$
and $\rho_0$ is invariant under (5.3). Define $v_0$ by (2.11).
The invariance means that $v_x(x,t)=\rho_0(x)$ a.e., and so
$$v(x,t)=\int_0^x \rho_0(y)dy+v(0,t)=v_0(x)+v(0,t)\,.
\tag 5.4
$$
Formula (5.3) operates like a semigroup, and this
with (5.4) implies that
$v(0,t)=Ct$ for some constant $C$. To determine 
$C$, we do a comparison. Let
$$v^e_0(x)=\rho^*x\ind\{x<0\}+(1-\rho^*)x\ind\{x>0\}\,,\,
v^e(x,t)=v^e_0(x)-t\lambda(0)/4
$$
and 
$$v^n_0(x)=(1-\rho^*)x\ind\{x<0\}+\rho^*x\ind\{x>0\}\,,\,
v^n(x,t)=v^n_0(x)-t\lambda(0)/4
$$
denote the evolution of the entropy and the non-entropy shock
at the origin. The bounds
$ v^n(x,t)\le v(x,t)\le v^e(x,t)$ are valid because they
are valid at time $t=0$  and preserved by (5.3). 
Consequently 
$$v(x,t)=v_0(x)-t\lambda(0)/4.
\tag 5.5
$$

\proclaim{Lemma 5.1} There cannot exist an $x\ne 0$,
and $\e,\delta>0$ such that this holds:
$$\aligned
&\text{$v_0(x)\ge v_0(q) +(\rho^*+\delta)(x-q)$ for
$q\in[x-\e,x]$}\\
\text{and }\,&\text{$v_0(x)\ge v_0(q) -(1-\rho^*-\delta)(q-x)$ for
$q\in[x,x+\e]$.}
\endaligned
\tag 5.6
$$
\endproclaim

\demo{Proof} Suppose such $x,\e,\delta$ exist. Pick $t>0$
small enough 
so that $t<\e$ and 
$Bt <|x|$. Then, by (5.2), (5.3) becomes
$$v(x,t)=\sup_{q\in\mmR}\bigl
\{v_0(q)-tg_0\bigl((x-q)/t\bigr)\bigr\}\,.
\tag 5.7
$$
We shall show that formula (5.7) gives something strictly
smaller than (5.5), and this contradiction makes (5.6)
impossible. 

In (5.7) it suffices to consider $q\in[x-t,x+t]$, by observing
from (2.13) that $g_0$ has constant slopes to the left of $-1$
and to the right of $1$. 
For $q\in[x,x+t]$ write 
$$\aligned
v_0(q)-tg_0\bigl((x-q)/t\bigr)
&\le v_0(x) + (1-\rho^*-\delta)(q-x)-tg_0\bigl((x-q)/t\bigr)\\
&=v_0(x) -t\bigl\{ (1-\rho^*-\delta)\xi + g_0(\xi)\bigr\}\,,
\endaligned
$$
where $\xi=(x-q)/t$ may vary freely in $[-1,0]$. 
By the duality of $g_0$ and the TASEP current $f_0(\rho)=\rho(1-\rho)$, 
the expression
in braces is bounded below by $f_0(1-\rho^*-\delta)$
$=$ $f_0(\rho^*+\delta)$. 
This and a similar argument for  $q\in[x-t,x]$ give 
$$v(x,t)\le v_0(x) -tf_0(\rho^*+\delta).$$
Since $\lambda(0)/4=f_0(\rho^*)$, the above bound is
 strictly less than (5.5), provided $\delta$ is
chosen small enough to have 
$\rho^*+\delta<1/2$. 
\qed
\enddemo

Now we can prove that $\rho_0$ takes only the values
$\{\rho^*, 1-\rho^*\}$, up to Lebesgue null sets. 
For suppose the set 
$A=\{ x: \rho^*+\delta_0 \le \rho_0(x)\le 1-\rho^*-\delta_0\}$
has positive Lebesgue measure for some $\delta_0>0$. 
Then $A$ has  a density point
$x$. (For a definition, see for example p.\ 107 in \cite{WZ}.)
We show that (5.6) holds at $x$. 
Let $q<x$. Since $\rho_0(x)\ge \rho^*$, 
$$\aligned
v_0(x)-v_0(q)&=\int_{[q,x]\cap A} \rho_0(y)dy
 + \int_{[q,x]\cap A^c} \rho_0(y)dy\\
&\ge \rho^*(x-q) +\delta_0\cdot m\bigl([q,x]\cap A\bigr)\,,
\endaligned
$$
where we wrote $m$ for Lebesgue measure.  
Since $x$ is a density point, $m\bigl([q,x]\cap A\bigr)\ge (x-q)/2$ if
$q$ is close enough to $x$. This checks the first part of
(5.6), and the other part is similar. 

Assuming that $\rho_0$ is piecewise continuous, we
now conclude that it is piecewise constant with 
values $\rho^*$ and $1-\rho^*$. 
To prove Cor.\ 2.2, it remains to observe that Lemma 5.1 
prevents a jump from $1-\rho^*$ to $\rho^*$ everywhere
else except at $x=0$. 

\hbox{}

\flushpar
{\it 5.3. Proof of Corollary 2.3 }

\hbox{}

Let the initial distribution of the process
be the i.i.d.\ product measure $\alpha_0$ with density
$\alpha_0\{\eta_i=1\}=\rho$, with $\rho$ outside the
disturbed range $(\rho^*,1-\rho^*)$. Let $\alpha_t$ be 
the distribution of the process at time $t$.
 Let $\mu$ be a 
limit point of the time averages of $\alpha_t$, so for
some sequence $t_k\nearrow \infty$, 
$$\mu=\lim_{k\to\infty} \frac1{t_k}\int_0^{t_k}\alpha_s ds$$
in the weak sense on the compact state space $\{0,1\}^\mmZ$. 
Such a limit point exists by compactness, and is then 
automatically invariant for the process. 

Write $E^{\alpha_0}$ for expectation under the path measure
of the process started with distribution $\alpha_0$. 
For $i\ne 0$ we have the equation
$$E^{\alpha_0}\bigl[J_i(t)\bigr]
=\int_0^t \alpha_s\{\eta_i=1,\eta_{i+1}=0\}ds
\tag 5.8
$$
 because $J_i(t)$ [$=$ the number
of jumps from site $i$ up to time $t$] increases by 1 at rate
1 when the event  inside the braces holds. Thus  
$$\aligned
\mu\{\eta_i=1,\eta_{i+1}=0\}&=
\lim_{k\to\infty} \frac1{t_k}\int_0^{t_k}\alpha_s 
\{\eta_i=1,\eta_{i+1}=0\}
ds\\
&=\lim_{k\to\infty} \frac1{t_k}E^{\alpha_0}\bigl[J_i(t_k)\bigr]\\
&=\rho(1-\rho)\,.
\endaligned
\tag 5.9
$$
The last equality is from (2.17) and Case 1 of Cor.\ 2.1, with initial
function $v_0(x)=\rho x$. 
$J_i(t)$ is bounded by a Poisson($t$) random variable, 
hence there is uniform integrability to justify the limit above. 
If $i=0$, we have to multiply
the right-hand side of (5.8)  by the factor
$r$, so the limit on the last
line of (5.9) is multiplied by $r^{-1}$. 
This proves (2.21) for $\mu$. 

Suppose the initial distribution $\alpha_0$ has a macroscopic
profile $\rho_0$ in the sense of (2.16) such that
$\rho^*\le \rho_0(x)\le 1-\rho^*$. Then the limit in
(5.9) is
valid with $\rho=\rho^*$. In particular, this is the 
case for the $\alpha_0$ with  non-entropy shock described in the 
second paragraph of Cor.\ 2.3.

\hbox{}

{\it Acknowledgements.} I thank Maury Bramson, Pablo Ferrari,
Tom Liggett,  and 
Fraydoun Rezakhanlou for discussions on this and related
problems, and  Dan Ostrov for guidance on Hamilton-Jacobi
equations with discontinuities. 

\hbox{}

\head References \endhead 

\hbox{}

\flushpar  
[CR]  Covert, P. and Rezakhanlou, F.: Hydrodynamic limit 
for particle systems with nonconstant speed parameter.
\jsp\ 88, 383--426 (1997).

\hbox{}

\flushpar  
[Ev] Evans, L. C.: Partial Differential Equations. American
Mathematical Society, 1998. 

\hbox{}

\flushpar  
[GW]  Glynn, P. W. and  Whitt, W.: 
Departures from many queues in a series.
 \aap\ 1, 546--572 (1991). 

\hbox{}

\flushpar  
[JL1] Janowsky, S. A. and Lebowitz, J. L.: Finite size
effects and shock fluctuations in 
the asymmetric simple exclusion process. 
Phys. Rev. A 45, 618--625 (1992). 

\hbox{}

\flushpar  
[JL2] Janowsky, S. A. and Lebowitz, J. L.: Exact results for
the asymmetric simple exclusion process with a blockage. 
\jsp\ 77, 35--51 (1994). 

\hbox{}

\flushpar  
[Jo] Johansson, K.: Shape fluctuations and random 
matrices. To appear in Comm. Math. Phys. 
Preprint math.CO/9903134. 

\hbox{}

\flushpar  
[Ki]  Kingman, J. F. C.:
The ergodic theory of subadditive
stochastic processes. 
J. Royal Stat. Soc. Ser. B 30  (1968) 499--510.

\hbox{}

\flushpar  
[KL] Kipnis, C. and  Landim, C.: Scaling Limits of
Interacting Particle Systems. 
Springer-Verlag,  New York 1999. 

\hbox{}

\flushpar  
[La]  Landim, C.: Hydrodynamical limit
for space inhomogeneous one-dimensional totally asymmetric
zero-range processes. \ap\  24 (1996)  599--638. 
   
\hbox{}

\flushpar  
[Li]  Liggett, T. M.: 
 Stochastic Interacting Systems: Contact, Voter
and Exclusion Processes. 
Springer-Verlag,  New York, 1999.

\hbox{}

\flushpar  
[Os] Ostrov, D.: Solutions of Hamilton-Jacobi 
equations and scalar conservation laws with discontinuous 
space-time dependence. Preprint (1999). 

\hbox{}

\flushpar  
[Ro] Rost, H.:  
Non-equilibrium behaviour of a many particle
process: Density profile and local equilibrium. 
\zwvg\  58, 41--53 (1981). 

\hbox{}

\flushpar  
[Se1] Sepp\"al\"ainen, T.:  Hydrodynamic scaling, convex duality,
and asymptotic shapes of growth models. 
 Markov Process.\ Related Fields.\ 4, 1--26 (1998).

\hbox{}

\flushpar  
[Se2] Sepp\"al\"ainen, T.: Coupling the totally asymmetric
simple exclusion process with a moving interface.
 Markov Process.\ Related Fields 4, 593--628 (1998). 

\hbox{}

\flushpar  
[Se3] Sepp\"al\"ainen, T.: Existence of hydrodynamics
for the totally asymmetric
simple $K$-exclusion process. 
 Ann.\ Probab. 27, 361--415 (1999).

\hbox{}

\flushpar  
[WZ] Wheeden, R. L. and Zygmund, A.: Measure and Integral. 
Marcel Dekker 1977.

\enddocument